\documentclass[a4paper,dvips,12pt]{amsart}

\usepackage{latexsym}
\usepackage{caption}
\usepackage[hang,tight,raggedright,TABTOPCAP,FIGTOPCAP]{subfigure}
\usepackage{paralist}
\usepackage[breaklinks=true,dvips]{hyperref}
\usepackage{psfrag}
\usepackage{graphicx}
\usepackage{hyperref}

\newtheorem{lemma}{Lemma}[section]

\newtheorem{proposition}[lemma]{Proposition}
\newtheorem{theorem}[lemma]{Theorem}
\newtheorem*{conjecture}{Conjecture}

\usepackage{txfonts}

\title{Lattice points in {M}inkowski sums}

\author[Ch.\ Haase]{Christian Haase}
\author[B.\ Nill]{Benjamin Nill}
\author[A.\ Paffenholz]{Andreas Paffenholz}
\author[F.\ Santos]{Francisco Santos}

\address{Inst. f\"ur Mathematik, Arnimallee 3, 14195 Berlin, Germany}
\email{\{christian.haase,nill,paffenho\}@math.fu-berlin.de}

\thanks{The first three authors were supported by Emmy Noether
  fellowship HA 4383/1 of the German research society DFG}

\address{Facultad de Ciencias, Universidad de Cantabria, Av. de los Castros s/n, 
E-39005 Santander, Spain}
\email{santosf@unican.es}

\thanks{The fourth author was supported by grant MTM2005-08618-C02-02 of the 
  Spanish Ministry of Science.}

\date{November 28, 2007}

\newcommand{\ov}[1]{\overline{#1}}
\numberwithin{equation}{section}

\newcommand{\clt}{\mathcal{L}}
\newcommand{\R}{\mathbb R}
\newcommand{\Z}{\mathbb Z}
\newcommand{\N}{\mathbb N}

\newcommand{\ms}{\mathfrak s}

\begin{document}

\begin{abstract}
  Fakhruddin has proved that for two lattice  polygons $P$ and $Q$ any
  lattice  point in their  Minkowski sum can be  written as a sum of a
  lattice point in $P$ and one in $Q$, provided $P$ is smooth and 
  the normal fan of $P$ is a subdivision of the normal fan of $Q$.

  We give a  shorter combinatorial proof  of this fact  that does not
  need the smoothness assumption on $P$.
\end{abstract}

\maketitle

\section{Introduction}
\noindent
It is one of those problems. Everyone can understand it
immediately. Yet, to this day we do not have any satisfactory
solution.

A lattice polygon $P \subset \R^2$ is the convex hull of finitely
many points in the lattice $\Z^2$.
Given two lattice polygons $P$ and $Q$, we consider the addition map
\begin{alignat*}{2}
  \ms:(P\cap \Z^2)\;\times\;& (Q\cap \Z^2)&&\;\longrightarrow\;
  (P+Q)\cap \Z^2 .\\
  (\,x,\;&y\,)&&\;\longmapsto\; x+y
\end{alignat*}
We want to understand when $\ms$ is surjective. Equivalently,
when is 
\[
(P\cap \Z^2)+(Q\cap \Z^2) =(P+Q)\cap \Z^2 \ ?
\]

This very basic question in discrete geometry (and its higher
dimensional analogue) appears in different guises in algebraic
geometry, commutative algebra, and integer programming. Specific
cases also arise in additive number theory, representation theory, and
statistics.
Motivation for a conjectured sufficient condition comes from 
algebraic geometry.
\begin{conjecture}[Oda]
  \label{p:problem1}
  Let $X$ be a smooth projective toric variety, let $D$ be an
  ample divisor on $X$, and let $D'$ be nef. Then, the following homomorphism
is surjective:
 \[
 H^0(X,\mathcal{O}(D)) \otimes H^0(X,\mathcal{O}(D')) \rightarrow
  H^0(X,\mathcal{O}(D+D')).
  \]
\end{conjecture}

The toric dictionary translates this into discrete
geometry as follows:

\begin{conjecture}[Oda'] \label{p:problem2}
  Let $P$ and $Q$ be lattice polytopes. If $P$ is smooth and
  the normal fan of $Q$ coarsens that of $P$, then the map $\ms$ is
  surjective.
\end{conjecture}

Here, a lattice polytope $P$ is called smooth if it is simple and
at every vertex the primitive facet normals generate the dual
lattice. This condition is equivalent to $P$ corresponding to an ample
divisor on a smooth toric variety. The case $Q=nP$, $n\in \N$ of this
conjecture is the conjecture that \emph{all smooth lattice polytopes
  are projectively normal}.

The two-dimensional case of Oda's conjecture is now Fakhruddin's
Theorem~\cite{0208.5178}, with an independent proof by
Ogata~\cite{1110.14050}.
In this note, we generalize Fakhruddin's Theorem to the non-smooth
case.
\begin{theorem}\label{thm:main}
  Let $P$ and $Q$ be lattice polygons such that the  normal fan of $Q$
  coarsens that of $P$.  Then  the map $\ms$ is surjective.
\end{theorem}

Our proof originated in a discussion about normality of polytopes
during a mini-workshop at Oberwolfach~\cite{owr_normality}.

\begin{figure}[ht]
  \centering
  \psfrag{P}[tc][tc]{$P$}
  \psfrag{Q}[tl][tl]{$Q$}
  \psfrag{mQ}[br][br]{$-Q$}
  \psfrag{PQ}[bl][bl]{$P+Q$}
  \includegraphics[height=.15\textheight]{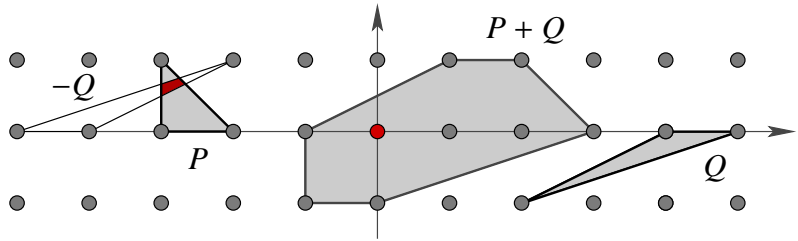}
  \caption{Theorem~\ref{thm:main} fails without  the assumption on the normal
    fan.}
  \label{fig:nfcond}
\end{figure}

The  assumption   on  the       normal fan   is    necessary.      See
Figure~\ref{fig:nfcond} for an example of two lattice polygons that do
not satisfy the condition on the normal fan.  The point $(0,0) \in
P+Q$ cannot be written as a sum of a lattice point in $P$ and one in
$Q$.

Embarrassingly, in dimension three the conjecture is open even if we take $Q=P$. 
Observe that in dimension three and higher the smoothness hypothesis cannot be removed. 
For example, if $P=Q$ is the simplex of lattice volume two in $\R^3$ having the four vertices 
as its only lattice points, 
then the centroid of $P+P$ is a lattice point but it is not in the image of the map $\ms$.

\section{Lattice point free intersections are $4$-gons}

As an intermediate step we prove the following curious result.

\begin{proposition}\label{prop:polygonintersection}
  Let $P$ and $Q$  be lattice  polygons and let $Z=P\cap Q$.
   If $Z$
  is not empty but does not contain a lattice point, then $Z$
  is a $4$-gon with two opposite edges coming from $P$ 
 and the other two coming from $Q$.
\end{proposition}
Examples of the stated $4$-gons  appear in Figures~\ref{fig:nfcond} and~\ref{fig:zlemma}.

\begin{proof}
  Let $Z := P\cap Q$.  If some vertex of $Z$ is a  vertex of $P$ or of
  $Q$, then it is a  lattice point in $Z$.  So, let us assume that all
  vertices of $Z$  arise   from an edge   of  $P$ and an  edge  of $Q$
  intersecting in their relative  interiors.  This implies that $Z$ is
  two-dimensional,  and that edges of  $Z$  are alternatingly edges of
  $P$ and of $Q$.  In particular, $Z$ has an even  number $n\geq 4$ of
  edges.  

  We prove the theorem  by contradiction.  For  this, assume $n\ge  6$
  and let $\clt(P)$ denote the  set of lattice points  in $P$ that are
  not  vertices of $P$.  We may  assume that $P$ minimizes $|\clt(P)|$
  among the polygons for which  $Z= P\cap Q$ has  more than four edges
  and contains no lattice point.

  If  $\clt(P)=\emptyset$, then $P$ is contained   in a (closed) strip
  $R$ of lattice width one.   The interior of $R$ intersects precisely
  two edges of  $Q$, since the strip  contains no lattice point in its
  interior.  Those two  are the only edges of  $Q$ that can contribute
  to edges  of  $Z$.    Hence,   $Z$ is   in  fact  a   $4$-gon.   See
  Figure~\ref{fig:zlemma}.

\begin{figure}[hb]
  \hfill
  \psfrag{Q}[bl][bl]{$P$}
  \psfrag{P}[bl][bl]{$Q$}
  \psfrag{R}[tr][tr]{$R$}
  \begin{minipage}[t]{.4\linewidth}
    \includegraphics[width=.85\linewidth]{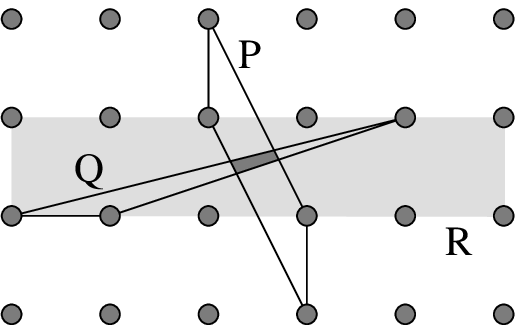}
    \caption{All lattice points of $P$ are vertices}
    \label{fig:zlemma}
  \end{minipage}
  \hfill
  \psfrag{q}[tr][tr]{$q$}
  % \psfrag{q2}{$q_2$}
  % \psfrag{q3}[tl][tl]{$q_3$}
  \psfrag{pl}[tl][tl]{$p_l$}
  \psfrag{pr}[br][br]{$p_r$}
  % \psfrag{p3}[tr][tr]{$p_3$}
  \psfrag{vl}[tl][tl]{$v_l$}
  % \psfrag{w1}[bl][bl]{$w_1$}
  \psfrag{vr}[b][b]{$v_r$}
  % \psfrag{w2}[tr][tr]{$w_2$}
  % \psfrag{v3}[tl][tl]{$v_3$}
  % \psfrag{w3}[tr][tr]{$w_3$}
  \psfrag{P}{$P$}
  \psfrag{Q}{$Q$}
  \psfrag{Z}{$Z$}
  \psfrag{m}{$m$}
  \psfrag{P1}[cl][cl]{$P^{(1)}$}
  \psfrag{Q1}[cl][cl]{$Q^{(2)}$}
  \begin{minipage}[t]{.45\linewidth}
    \centering
    \includegraphics[width=.9\linewidth]{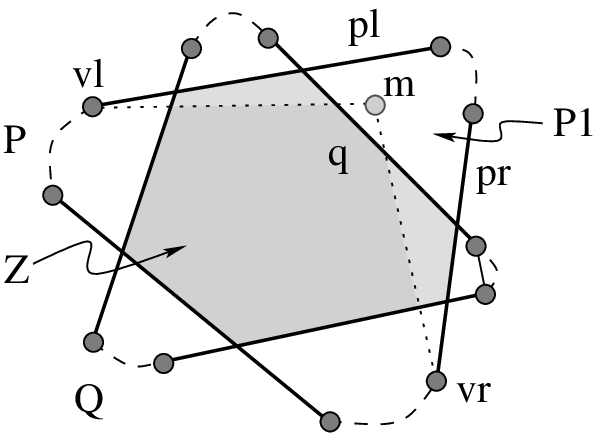}
    \caption{$m\in P\cap \Z^2$ is not a vertex of $P$}
    \label{fig:zgen}
  \end{minipage}
  \hfill
\end{figure}

Now assume $|\clt(P)|>0$.  We will construct a subpolytope $P' \subset
P$   with $|\clt(P')|  < |\clt(P)|$ and   such  that the  intersection
$Z'=P'\cap Q$ has the same number of edges as $Z$, a contradiction.

For this, let $m\in\clt(P)$. 
By assumption, $m\not \in Z$.
Hence, there is an edge $q$ of $Z$ with $m$ in its outer half-space
$H^+$. As $q$ comes from an edge of $Q$, both $Q$ and $Z$ are
contained in the closed half-space $H^-$. See Figure~\ref{fig:zgen}.

  Let  $p_l, p_r$  be the  edges of $Z$   adjacent to  $q$. Then $p_l$
  (respectively, $p_r$) is part of  an edge $\ov{p_l}$  (respectively,
  $\ov{p_r}$) of $P$. Let $v_l$ (respectively, $v_r$) be the vertex of
  $\ov{p_l}$ (respectively, $\ov{p_r}$) contained in $H^-$.  Since $Z$
  is not a $4$-gon we have $v_l \not= v_r$.

  We  define $P'$ as the convex  hull of $m$ and   all vertices of $P$
  that are contained  in $H^-$.  By construction,  $Z':=P'\cap Q$ is a
  $2$-dimensional polygon with the same number of vertices as $Z$.  As
  $m$ is a vertex of $P'$, $\clt(P')\subsetneq \clt(P)$
  \end{proof}

\section{Proof of Theorem~\ref{thm:main}}

We first translate Theorem~\ref{thm:main} into a  statement that
does not  involve  the map  $\ms$  anymore.  The following 
  necessary and
sufficient condition for $\ms$ to be surjective is due to
 Benjamin Howard~\cite{benHoward}.
 
\begin{lemma}
  $\ms$ is surjective if and only if for all $z\in\Z^2$
  \begin{align*}
    P\cap  (z-Q)\ne\emptyset\quad\Longleftrightarrow\quad  \left(P\cap
      (z-Q)\right)\cap\Z^2\ne\emptyset
  \end{align*}
\end{lemma}

 \begin{proof}
   The left-hand side is equivalent to $z\in P+Q$. The right-hand side
   to $z\in (P\cap \Z^2) +(Q\cap \Z^2)$.
 \end{proof}

For example, in Figure~\ref{fig:nfcond} the point $(0,0)$ is not  in
the image of $\ms$ because the intersection of $P$  and $-Q$ does not
contain a lattice point.  
Using this lemma,  Theorem~\ref{thm:main} is equivalent to the following.

\begin{theorem}\label{thm:reformulated}
  Let $P$ and $Q$ be lattice polygons. If $P\cap Q \not= \emptyset$ and the normal fan of $-Q$
  coarsens that of $P$, then $P \cap Q$ 
  contains a  lattice point.
\end{theorem}

\begin{figure}[ht]
  \begin{minipage}[t]{.48\linewidth}
    \centering \psfrag{Z}{$Z$} \psfrag{e1}[bl][bl]{$e_1$}
    \psfrag{e2}[tr][tr]{$e_2$} \psfrag{f1}[tl][tl]{$f_1$}
    \psfrag{f2}[rc][rc]{$f_2$}
    \psfrag{of2}[br][br]{$\ov f_2$}
    \psfrag{of1}[br][br]{$\ov f_1$}
    \includegraphics[width=.8\textwidth]{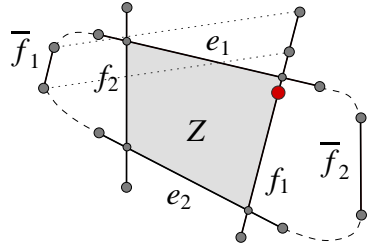}
    \caption{A lattice point on $f_1$.}
    \label{fig:lp_on_edge}
  \end{minipage}
\end{figure}

\begin{proof}
Let $Z=P\cap Q$ and suppose that $Z\cap \Z^2$ was empty. Then, by 
Proposition~\ref{prop:polygonintersection},
$Z$ is a $4$-gon with two opposite edges coming from $P$ and
the other two coming from $Q$.

  Let $e_1$ and $e_2$ be the edges of $Z$ originating from $P$ and $f_1,
  f_2$ those from $Q$.  Let $\ov f_1$  and $\ov f_2$ be the edges of $P$ with 
  exterior normals opposite to those of $f_1$ and $f_2$, respectively.
  They exist by the condition on the  normal fans.
  
  By construction, $\ov f_2$, $f_1$, $f_2$ and $\ov f_1$ are all
  contained, and appear in this order, in the (possibly degenerate)
  wedge defined by  the lines  supporting $e_1$ and  $e_2$.
  See Figure~\ref{fig:lp_on_edge}.
  In particular, at least one of the $\ov f_i$ is shorter than or
  equal to its corresponding $f_i$. But since $\ov f_i$ is a lattice
  segment, and since $f_i$ is parallel to it, contained in a lattice
  line, and equal to or longer than it, $f_i$ must contain a lattice
  point. This contradiction finishes the proof.
\end{proof}
  
\bibliographystyle{alpha}

\begin{thebibliography}{HHM07}

\bibitem[Fak02]{0208.5178}
Najmuddin Fakhruddin.
\newblock {Multiplication maps of linear systems on smooth projective toric
  surfaces.}
\newblock{Preprint, \href{http://arxiv.org/abs/math.AG/0208178}{math.AG/0208178}}, 2002.

\bibitem[HHM07]{owr_normality}
Christian Haase, Takayuki Hibi, and Diane MacLagan, editors.
\newblock {\em Mini-Workshop: Projective normality of smooth toric varieties},
  volume~39 of {\em Oberwolfach report}, 2007.

\bibitem[How07]{benHoward}
Benjamin~J. Howard.
\newblock {Matroids and geometric invariant theory of torus actions on flag
  spaces.}
\newblock {\em J. Algebra}, 312(1):527--541, 2007.

\bibitem[Oga06]{1110.14050}
Shoetsu Ogata.
\newblock {Multiplication maps of complete linear systems on projective toric
  surfaces.}
\newblock {\em Interdiscip. Inf. Sci.}, 12(2):93--107, 2006.

\end{thebibliography}

\setlength{\parindent}{0pt}

\end{document}